\documentstyle[12pt]{amsart}

\makeatletter
\def\cases#1{\left\{\,\vcenter{\normalbaselines\m@th
    \ialign{$##\hfil$&\quad##\hfil\crcr#1\crcr}}\right.}
\makeatother

\author{Jie Wu}
\address{MSRI\\
1000 Centennial Drive\\
Berkeley, CA 94720-5070}
\title{On Combinatorial Descriptions of Homotopy Groups of $\Sigma K(\pi,1)$}
\thanks{Research at MSRI is supported in part by NSF grant DMS-9022140}

\newtheorem{theorem}{Theorem}[section]
\newtheorem{definition}[theorem]{Definition}
\newtheorem{definitions}[theorem]{Definitions}

\newtheorem{notations}[theorem]{Notations}

\newtheorem{example}[theorem]{Example}

\newtheorem{lemma}[theorem]{Lemma}
\newtheorem{remark}[theorem]{Remark}

\newtheorem{proposition}[theorem]{Proposition}
\newtheorem{corollary}[theorem]{Corollary}

\begin{document}
\maketitle
\begin{abstract}
We give a combinatorial description of homotopy groups of $\Sigma K(\pi,1)$. In particular, all of the homotopy groups of the $3$-sphere are combinatorially given.
\end{abstract}
\section{Introduction}

In this article, we study homotopy groups and some related problems by using simplicial
 homotopy theory. The point of view here is that combinatorial aspects of group theory provide
further information about homotopy groups.\\
\par
The homotopy groups of the $3$-sphere, the suspension of $K(\pi,1)$ and
wedges of $2$-spheres are shown to be the centers of certain groups with specific generators 
and specific relations.  We list two group theoretical descriptions of $\pi_*(S^3)$ as follows.

\begin{definition}
{\bf A bracket arrangement of weight $n$} in a group is a set of elements defined recursively
as follows:
\par
 Let $G$ be a group and let $a_1,a_2,\cdots,a_n$ be a finite sequence of elements of $G$.
Let
$$
\beta^1(a_1)=a_1
$$
and if $n>1$, then let
$$
\beta^n(a_1,\cdots,a_n)=[\beta^k(a_1,\cdots,a_k),\beta^l(a_{k+1},\cdots,a_n)]
$$
for some $k$ and $l>0$.
\end{definition}
\begin{theorem}
For $n\geq 1$, $\pi_{n+2}(S^{3})$ is isomorphic to the center of the group with generators
 $y_{0},\dots,y_{n}$ and relations 
$$[y_{i_{1}}^{\epsilon_{1}},y_{i_{2}}^{\epsilon_{2}},\dots,y_{i_{t}}^{\epsilon_{t}}]$$
 with $\{i_{1},\dots,i_{t}\}=\{-1,0,\dots,n\}$ as sets in which the
 indices $i_{j}$ can be repeated, where 
$\epsilon_{j}=\pm 1, y_{-1}=(y_{0}\dots y_{n})^{-1}$ and 
the commutator bracket $[ \dots ]$ runs over all bracket arrangements of weight $t$ for
 each $t$.
\end{theorem}

\begin{notations}
Let $G$ be a group and let $S$ be a subset of $G$. Let $<\!S\!>$ denote the normal subgroup generated by
 $S$. Let $H_{j}$ be a sequence of subgroups of $G$ for $1\leq j\leq k$. Let
$[[H_{1},\dots,H_{k}]]$ denote the subgroup of $G$ generated by all of 
the commutators $[h_{i_1}^{(1)},\dots,h_{i_t}^{(t)}]$ with 
$\{i_1,\dots,i_t\}=\{1,\dots,k\}$ as sets in which the indices $i_j$can be repeated, 
where $h_j^{(s)}\in H_j$ and the commutator bracket
 $[ \dots ]$ runs over all of the bracket arrangements of weight $t$ for each $t$. 
\end{notations}

\begin{theorem}
Let $n\geq 1$. In the free group $F(y_0,\dots,y_n)$ freely generated
by $y_0,\allowbreak \dots,\allowbreak y_n$ we have
 $$
 ([[<\!y_0\!>,<\!y_1\!>,\dots,<\!y_n\!>]]\cap <\!y_{-1}\!>)/ [[<\!y_{-1}\!>,<\!y_0\!>,\dots,<\!y_n\!>]]\cong \pi_{n+2}(S^{3})
 $$
where $y_{-1}=(y_0\dots y_n)^{-1}$.
\end{theorem}

The method of the proofs of these theorems is to study the Moore chain complex of Milnor's
construction $F(S^1)$ for the $1$-sphere $S^1$.  A group theoretical description of the homotopy groups $\pi_*(\Sigma K(\pi,1))$
is as follows.

\begin{theorem}
Let $\pi$ be any group and let $\{x^{(\alpha)}|\alpha\in J\}$ be a set of generators for $\pi$. Then, for
$n\not=1$, $\pi_{n+2}(\Sigma K(\pi,1))$ is isomorphic to the center of the quotient group
of the free product 
$${\coprod^{groups}}_{0\leq j\leq n}(\pi)_j$$
modulo the relations
$$
[y^{(\alpha_1)\epsilon_1}_{i_1},y^{(\alpha_2)\epsilon_2}_{i_2},\cdots,y^{(\alpha_t)\epsilon_t}_{i_t}]
$$
with $\{i_1,i_2,\cdots,i_t\}=\{-1,0,1,\cdots,n\}$ as sets,
where $(\pi)_j$ is a copy of $\pi$ with generators $\{x^{(\alpha)}_j|\alpha\in J\}$, $\epsilon_j=\pm1$,
$y^{(\alpha)}_{-1}={x^{(\alpha)}_0}^{-1}$, $y_j=x^{(\alpha)}_j{x^{(\alpha)}_{j+1}}^{-1}$ for $1\leq j\leq n-1$, $y^{(\alpha)}_n=x^{(\alpha)}_m$ and the commutator bracket
$[\cdots]$ runs over all of the commutator bracket arrangment of weight $t$ for each $t$.
\end{theorem}

In particular, for $\pi={\bf Z}/m$, we have
\begin{corollary}
$\pi_{n+2}(\Sigma K({\bf Z}/m,1))$ is isomorphic to the center of the group with generators
$x_{0},\dots,x_{n}$
 and relations:
(1).
$x_j^m$
for $0\leq j\leq n$ and
(2). 
$$[y_{i_{1}}^{\epsilon_{1}},y_{i_{2}}^{\epsilon_{2}},\dots,y_{i_{t}}^{\epsilon_{t}}]$$
 with $\{i_{1},\dots,i_{t}\}=\{-1,0,\dots,n\}$ as sets in which the
 indices $i_{j}$ can be repeated, where 
$\epsilon_{j}=\pm 1$, $y_{-1}=x_0^{-1}$, $y_j=x_{j-1}x_j^{-1}$ for $-1\leq j\leq n-1$, $y_n=x_n$ and 
the commutator bracket $[ \dots ]$ runs over all bracket arrangements of weight $t$ for
 each $t$.
\end{corollary}

One of the features of Theorems 1.2 and 1.5 is that homotopy groups embed in certain 
'enveloping groups'. These 'enveloping groups' have systematic and uniform structure. The centers of 
these groups are of course more complicated to analyze. Theorem 1.4 is very similar to the
combinatorial decription of J. H. C. Whitehead's conjecture [see, Bo2, pp.317]. These decriptions give a
combinatorial question how to give a computable way to understand the quotient groups
$R\cap S/[R,S]$ for certain subgroups $R$ and $S$ of a finite generated free group. Few informations
about this question are known [see, Bo1]. The article is organized as follows.\\
\par
In Section 2, we study some general properties of simplicial groups. Central extensions in the
Moore-Postnikov systems will be considered. In Section 3, we study the intersection of certain
subgroups in free groups. The proofs of Theorems 1.2 and 1.4 is given in Section 4, where
Theorem 4.5 is Theorem 1.2 and Theorem 4.10 is Theorem 1.4. The proof
of Theorem 1.5 is given in Section 5, where Theorem 5.9 is Theorem 1.5.
In Section 6, we give some applications of our descriptions. One example is to compute
the homotopy groups of the Cohen's construction on the $1$-sphere. A direct corollary is to
give a short proof of the Milnor's counter example for the minimal simplicial groups.\\
\par
The author would like to thank F. Cohen, B. Gray, J. Harper, D. Kan, M. Mahowald, J. Moore, P. May, S. Priddy 
and many other mathematicians for their kindly encouragement and helpful suggestions.
 The author is indebted to helpful discussions with J. C. Moore and F. R. Cohen.

\section{Central Extensions in simplicial group theory}

In this section, we study some general properties of simplicial groups.
 A simplicial set $K$ is called a Kan complex if it satisfies 
the ${\it extension\, condition}$,i.e, for any simplicial map $f:\Lambda^k[n]\to K$ has
 an extension $g:\Delta[n]\to K$, where $\Delta[n]$ is the standard $n$ simplex, 
$\Lambda^k[n]$ is the subcomplex of $\Delta[n]$ generated by all
 $d_i(\sigma_n)$ for $i\not=k$ and $\sigma_n$ is the nondegenerate $n$ simplex 
in $\Delta[n]$,[K1, C2]. Recall that any simplicial group is a Kan complex [Mo1].
Given a simplicial group $G$, the Moore chain complex $(NG,d_0)$ is defined by
 $NG_n=\bigcap_{j\not=0}Kerd_j$ with $d_0:NG_n\to NG_{n-1}$. The classical Moore Theorem is
that $\pi_n(G)\cong H_n(NG)$ [Mo1, Theorem 4;C2, Theorem 3.7  or K2, Proposition 5.4]. 
Now let ${\cal Z}_n={\cal Z}_n(G)=\bigcap_jKerd_j$ denote the cycles and let
 ${\cal B}_n={\cal B}_n(G)=d_0(NG_{n+1})$ denote the boundaries.
It is easy to check that ${\cal B}_n$ is a normal subgroup of $G_n$ for any simplicial group $G$.

\begin{lemma}
Let $G$ be a  simplicial group. Then the homotopy group $\pi_n(G)$ is contained in the center
of $G_n/{\cal B}G_n$ for $n\geq1$.
\end{lemma}

\noindent{\em Proof.} Notice that $\pi_n(G)\cong{\cal Z}G_n/{\cal B}G_n$. It suffices to show that
the commutator $[x,y]\in{\cal B}G_n$ for any $x\in{\cal Z}G_n$ and $y\in G_n$.
\par
Notice that $x$ is a cycle. There is a simplicial map $f_x: S^n\rightarrow G$ such that
$f_x(\sigma_n)=x$, where $S^n$ is the standard $n$-sphere with a nondegenerate $n$-simplex $\sigma_n$.
 Now let the simplicial map $f_y:\Delta[n]\rightarrow G$ be the representative of $y$,
i.e, $f_y(\tau_n)=y$ for the nondegenerate $n$-simplex $\tau_n$. Let $\phi$ be the composite
$$
S^n\stackrel{j}{\rightarrow} S^n\wedge\Delta[n]\stackrel{[f_x,f_y]}{\rightarrow}G,
$$
where $j(\sigma_n)=\sigma_n\wedge\tau_n$ and $[f_x,f_y](a\wedge b)=[f_x(a),f_y(b)]$, the commutator
of $f_x(a)$ and $f_y(b)$. Notice that $\phi(\sigma_n)=[x,y]$ and $S^n\wedge\Delta[n]$ is contactible.
Thus $[x,y]\in{\cal B}_n$. The assertion follows.

\begin{definitions}
Let $G$ be a simplicial group. The subsimplicial group {\bf $R_nG$} is defined by setting
$$
(R_nG)_q=\{x\in G_q|f_x(\Delta[q])^{[n]}=1\},
$$
where $f_x$ is the representative of $x$ and $X^{[n]}$ is the $n$-skeleton of the simplicial set $X$.
The subsimplicial group {\bf $\overline{R}_nG$} is defined by setting
$$
(\overline{R}_nG)_q=\{x\in G_q|f_x((\Delta[q])_n)\subseteq{\cal B}G_n\}.
$$
Let {\bf $P_nG$} denote $G/R_nG$ and let {\bf $\overline{P}_nG$} denote $G/\overline{R}_nG$.
\end{definitions}

It is easy to check that both $R_nG$ and $\overline{R}_nG$ are normal subsimplicial group of $G$. Thus
both $P_nG$ and $\overline{P}_nG$ are quotient simplicial group of $G$. Notice that
$$
R_nG\subseteq\overline{R}_nG\subseteq R_{n-1}G.
$$ 
The quotient simplicial group $\overline{P}_nG$ is between $P_nG$ and $P_{n-1}G$. By checking the definition
of Moore-Post\-nikov systems of a simplicial set [Mo1;C2], we have
\begin{lemma}
The quotient simplicial group $P_nG$ is the standard $n$-th Moore-Postnikov system of the simplicial group
$G$.
\end{lemma}
The quotient simplicial group $\overline{P}_nG$ has the same homotopy type of $P_nG$.
\begin{proposition}
The quotient simplicial homomorphism $q_n:P_nG\rightarrow\overline{P}_nG$ is a homotopy equivalence for each $n$.
\end{proposition}
\noindent{\em Proof:}
The Moore chain complex of $P_nG$ is as follows:
$$
N(P_nG)_q=\cases{1&for\,$q>n+1$,\cr NG_{n+1}/{\cal Z}G_{n+1}&for\,$q=n+1$,\cr NG_q&for\,$q<n+1$.\cr}
$$
The Moore chain complex of $\overline{P}_nG$ is as follows:
$$
N(\overline{P}_nG)_q=\cases{1&for\,$q>n$,\cr NG_n/{\cal B}G_n &for\, $q=n$,\cr NG_q &for $q<n$.\cr}
$$ 
Thus $(q_n)_*:\pi_*(P_nG)\rightarrow\pi_*(\overline{P}_nG)$ is an isomorphism and the assertion follows.\\
\par
Let $\overline{F}_n$ denote the kernel of the quotient simplicial homomorphism 
$r_n:\overline{P}_nG\rightarrow P_{n-1}G$ for $n>0$.
\begin{proposition}
The simplicial group $\overline{F}_n$ is the minimal simplicial group $K(\pi_nG,n)$ for each $n>0$.
\end{proposition}
\noindent{\em Proof:} It is directly to check that
$$
N(\overline{F}_nG)_q=\cases{1&for\,$q\not=n$,\cr \pi_nG&for\, $q=n$.\cr}
$$
The assertion follows.

\begin{proposition}
 The short exact sequence of simplicial groups
$$
0\rightarrow\overline{F}_nG\rightarrow\overline{P}_nG\rightarrow P_{n-1}G\rightarrow0
$$
is a central extension for each $n>0$.
\end{proposition}
\noindent{\em Proof:} Consider the relative commutator sussimplicial group $[\overline{F}_nG,\overline{P}_nG]$.
By Lemma 2.1, $[\overline{F}_nG,\overline{P}_nG]_n=1$. Notice that $[\overline{F}_nG,\overline{P}_nG]$
is a subsimplicial group of minimal simplicial group $\overline{F}_nG\cong K(\pi_nG,n)$. Thus
$[\overline{F}_nG,\overline{P}_nG]=1$ and the assertion follows.
 
\begin{definition}
A simplicial group is said to be {\it $r$-centerless} if  the center $Z(G_n)=\{1\}$ for $n\geq r$.
\end{definition}

\begin{proposition}
Let $G$ be a reduced $r$-centerless simplicial group. Then $\pi_n(G)\cong Z(G_n/{\cal B}_n)$
for $n\geq r+1$
\end{proposition}
\noindent{\em Proof.} By Lemma 2.1, ${\cal Z}_n/{\cal B}_n\subseteq Z(G_n/{\cal B}_n)$.
It suffices to show that $Z(G_n/{\cal B}_n)\subseteq{\cal Z}_n/{\cal B}_n$ for $n\geq r+1$.
Now let ${\tilde{x}}\in Z(G_n/{\cal B}_n)$. Choose $x\in G_n$ with $p(x)={\tilde x}$, where
$p:G_n\to G_n/{\cal B}_n$ is the quotient homomorphism. To check
 ${\tilde x}\in{\cal Z}_n/{\cal B}_n$, it suffices to show that $x\in {\cal Z}_n$ or $d_jx=1$ for all $j$.
Since $Z(G_{n-1})=\{1\}$, $d_jx=1$ if and only if $[d_jx,y]=1$ for all $y\in G_{n-1}$. Now
$[d_jx,y]=d_j[x,s_{j-1}y]$ for $j>0$ and $[d_0x,y]=d_0[x,s_0y]$. Since
 ${\tilde x}\in Z(G_n/{\cal B}_n)$, $[x,z]\in{\cal B}_n\subseteq{\cal Z}_n$ for all $z\in G_n$
and therefore $[d_jx,y]=1$ for all $y\in G_{n-1}$. The assertion follows.\\
\par
By inspecting the proof, we also have

\begin{proposition}
Let $G$ be a reduced $r$-centerless simplicial group. Then 
$$Z(G_n/{\cal Z}_n)=\{1\}$$ for $n\geq r+1$.
\end{proposition}

\begin{lemma}
Let $G$ be a reduced simplicial group so that $G_n$ is cyclic or centerless for each n. Then
there exists a unique integer $\gamma_G>0$ so that $G_n=\{1\}$ for $n<\gamma_G$ and $Z(G_n)=\{1\}$
for $n>\gamma_G$.
\end{lemma}

\noindent{\em Proof.} Let $\gamma_G=max\{\gamma|G_n=\{1\}\,\,for\,\,n<\gamma\}$. Then $\gamma_G>0$. If
$\gamma_G<\infty$. Then $G_{\gamma_G}\neq\{1\}$. We show that $G_{\gamma_G+q}$ is centerless
for each $q>0$. Notice that 
$d_0^q\circ s_0^q:G_n\rightarrow G_{n+q}\rightarrow G_n$ and 
$d_1^q\circ s_1^q:G_n\rightarrow G_{n+q}\rightarrow G_n$ are identities for $n>0$
Thus $s_0^q(G_{\gamma_G})$ and $s_1^q(G_{\gamma_G})$ are nontrivial summands of 
$G_{\gamma_G+q}$. Now let $x=s_0^qy=s_1^qz\in s_0^q(G_{\gamma_G})\cap s_1^q(G_{\gamma_G})$.
Then $d_{q+1}x=d_{q+1}s_0^qy=s_0^qd_1y=1=d_{q+1}s_1^qz=s_1^{q-1}d_2s_1z=s_1^{q-1}z$.
Thus $x=s_1^qz=1$. And therefore $s_0^q(G_{\gamma_G})\cap s_1^q(G_{\gamma_G})=\{1\}$.
The assertion follows.

\begin{corollary}
Let $G$ be a reduced simplicial group such that $G_n$ is cyclic or centerless for each $n$.
Then $\pi_n(G)\cong Z(G_n/{\cal B}_n)$ for $n\neq\gamma_G+1$, where $\gamma_G$ is defined as above.
\end{corollary}

Notice that, for any free group $F$, $rank(F)\geq2\Leftrightarrow Z(F)=\{1\}$ and
$F\neq\{1\}$. We have
\begin{lemma}
Let $G$ be a reduced simplicial group such that $G_n$ is a free group for each $n$. Then there
exits a unique integer $\gamma_G>0$ so that $G_n=\{1\}$ for $n<\gamma_G$ and $rank(G_n)\geq2$
for $n>\gamma_G$.
\end{lemma}

\begin{example}
The $1$-stem is determined in this example.
\end{example}
 Let $G=F(S^n)$,  Milnor's $F$-construction on the standard n-sphere for $n\geq1$.
Then $G_n\cong F(\sigma)\cong {\bf Z}(\sigma)$, the free abelian group
generated by $\sigma$, 
$$G_{n+1}\cong F(s_0\sigma,s_1\sigma,\dots,s_n\sigma)$$
and $G_{n+2}\cong F(s_is_j\sigma |0\leq j<i\leq n)$.
 It is easy to check that $\Gamma^2G_{n+1}={\cal Z}_{n+1}$,
where $\Gamma^qG$ is the $q$-th term in the lower central series
of a group $G$ starting with $\Gamma^1G=G$. 
By Lemma 2.1,
 $$\Gamma^3G_{n+1}=[{\cal Z}_{n+1}, G_{n+1}]\subseteq{\cal B}_{n+1}.$$

If $n=1$, then it will be shown that $\Gamma^3G_{n+1}={\cal B}_{n+1}$ in Section 4 and therefore
 $\pi_3(S^2)\cong\pi_2(F(S^1))\cong{\bf Z}$, which is generated by $[s_0\sigma,s_1\sigma]$.
\par
Suppose that $n>1$. Consider the following equations 
$$
d_k([s_{j-1}s_i\sigma,s_{j+1}s_j\sigma])=\cases{1&for\,$k\neq j,$\cr[s_i\sigma,s_j\sigma]&for\,k=j,\cr}
$$
for $i+1<j\leq n$,
$$
d_k[s_{i+2}s_{i+1}\sigma,s_{i+3}s_i\sigma]=\cases{[s_{i+1}\sigma,s_{i+2}\sigma]&for\,k=i+1,\cr[s_{i+1}\sigma,s_i\sigma]&for\,k=i+3,\cr1&otherwise,\cr}
$$
and
$$
d_k[s_{i+2}s_i\sigma,s_{i+3}s_{i+1}\sigma]=\cases{[s_{i+1}\sigma,s_{i+2}\sigma]&for\,k=i+1,\cr[s_i\sigma,s_{i+2}\sigma]&for\,k=i+2,\cr[s_i\sigma,s_{i+1}\sigma]&for\,k=i+3,\cr1&otherwise.\cr}
$$
By the Homotopy Addition Theorem [C2, Theorem 2.4], $[s_i\sigma, s_j\sigma]\in{\cal B}_{n+1}$ for $i+1<j$,
 $[s_{i+1}\sigma,s_{i+2}\sigma]\equiv[s_i\sigma,s_{i+1}\sigma]$ mod ${\cal B}_{n+1}$ and
$0\equiv[s_i\sigma,s_{i+2}\sigma]\equiv[s_{i+1}\sigma,s_{i+2}\sigma]\equiv 2[s_i\sigma,s_{i+1}\sigma]$  
if $i+2\leq n$. Notice that $[s_0\sigma,s_1\sigma]\neq0$ in $\pi_{n+1}(\Gamma^2G/\Gamma^3G)$.
Thus $[s_0\sigma,s_1\sigma]\notin{\cal B}_{n+1}$ and, by the relations above,
$$
\pi_{n+2}(S^{n+1})\cong\pi_{n+1}(G)\cong{\cal Z}_{n+1}/{\cal B}_{n+1}\cong{\bf Z}/2
$$
 for $n\geq2$,
which can be represented by $[s_i\sigma,s_{i+1}\sigma]$ for $0\leq i\leq n-1$.

\section{Intersections of certain subgroups in free groups}
In this section, we give some group theory preliminary. The intersections of certain subgroups in free
groups are considered in this section. We will use these informations to determine the Moore chain complex
of certain simplicial groups in the next sections.
\begin{definition}
 let $S$ be  a set and let $T\subseteq S$ a subset. The {\bf projection homomorphism}
 $$\pi:F(S)\rightarrow F(T)$$
is defined by 
$$
\pi(x)=\cases{x&$x\in T$,\cr 1&$x\notin T$.\cr}
$$  
\end{definition}

Now let $\pi:F(S)\rightarrow F(T)$ be a projection homomorphism and let $R$ equal
 the kernel of $\pi$. Define the subsets of the free group $F(S)$ as follows.
$$
{\cal A}_T(k)=\{[[x,y_1^{\epsilon_1}]\cdots],y_t^{\epsilon_t}]|0\leq t\leq k,\epsilon_j=\pm1,
y=y_1^{\epsilon_1}\dots y_t^{\epsilon_t}\in F(T),x\in S-T\},
$$
where $y=y_1^{\epsilon_1}\cdots y_t^{\epsilon_t}\in F(T)$ runs over reduced words in $F(T)$ with $t\leq k$ and
$y_j\in T$. Furthermore define $[[x,y_1^{\epsilon_1}]\dots],y_t^{\epsilon_t}]=x$ for $t=0$.
Define 
$$
{\cal B}_T(k)=\{\phi^{-1}x\phi|\phi\in F(T)\,a\,reduced\,word\,with\,lenth\,l(\phi)\leq k,x\in S-T\},
$$

$${\cal A}_T=\cup_{k\geq0}{\cal A}_T(k)$$
 and 
$${\cal B}_T=\cup_{k\geq0}{\cal B}_T(k).$$
 By the classical Kurosch-Schreier
theorem ( see [MKS, pp.243, K2, Theorem 18.1]), we have

\begin{proposition}
The subgroup $R$ is a free group freely generated by ${\cal B}_T$.
\end{proposition}

We will show that ${\cal A}_T$ is also a set of free generators for $R$. We need a lemma. 

\begin{lemma}
Let $\phi:F_1\rightarrow F_2$ be a  homomorphism of free groups. Suppose that $\phi^{ab}:F^{ab}_1\rightarrow F^{ab}_2$
is an isomorphism, where $F^{ab}$ is the abelianlizer of the group $F$. Then $\phi:F_1\rightarrow F_2$ is a monomorphism. 
\end{lemma}

\noindent{\em Proof.} Notice that $\phi_*: H_*(F_1)\rightarrow H_*(F_2)$ is an isomorphism,
where $H_*(G)$ is the homology of the group $G$. Thus
$F_1/\Gamma^rF_1\rightarrow F_2/\Gamma^r F_2$ is an isomorphism for each $r$, where $\Gamma^rG$ is the $r$-th
term in the lower central seris of the group $G$ starting with $\Gamma^1G=G$ and so
$$ lim_r F_1/\Gamma^r F_1\rightarrow lim_r F_2/\Gamma^r F_2$$
is an isomorphism. Notice that $\cap_r\Gamma^rF=1$ for any free group $F$. Thus $F\rightarrow lim_r F/\Gamma^r F$ is a monomorphism. The assertion follows.\\
\par
\begin{proposition}
The subgroup $R$ is a free group freely generated by ${\cal A}_T$.
\end{proposition}

\noindent{\em Proof.} First we assume that both $S$ and $T$ are finite sets. Denote by 
$i_k:{\cal A}_T(k)\rightarrow R$ and $j_k:{\cal B}_T(k)\rightarrow R$ the natural inclusions. Notice that
$$R=F({\cal B}_T)=colim_kF({\cal B}_T(k)).$$ We set up the following steps.\\
\par
\noindent{\em Step1.} ${\cal A}_T(k)\subseteq F({\cal B}_T(k))$.\\
\par
The proof of this statement is given by induction on $k$ starting with ${\cal A}_T(0)={\cal B}_T(0)=S-T$. Suppose that 
${\cal A}_T(k-1)\subseteq F({\cal B}_T(k-1))$ and let
 $w=[[x,y_1^{\epsilon_1},\cdots,y_t^{\epsilon_t}]\in{\cal A}_T(k)$. If $t<k$, then
 $w\in{\cal A}_T(k-1)\subseteq F({\cal B}_T(k-1))\subseteq F({\cal B}_T(k))$, by induction.
Now 
$$
[[x,y_1^{\epsilon_1}]\cdots],y_k^{\epsilon_k}]=
[[x,y_1^{\epsilon_1}]\cdots],y_{k-1}^{\epsilon_{k-1}}]^{-1}\cdot y_k^{-\epsilon_k}[[x,y_1^{\epsilon_1}]\cdots],y_{k-1}^{\epsilon_{k-1}}]y_k^{\epsilon_k}
$$
Now since $[[x,y_1^{\epsilon_1}]\cdots],y_{k-1}^{\epsilon_{k-1}}]\in F({\cal B}_T(k-1))$,
$
[[x,y_1^{\epsilon_1}]\cdots],y_{k-1}^{\epsilon_{k-1}}]=\prod_{j=1}^s(\phi_j^{-1}x_j\phi_j))^{\eta_j}
$
with $\phi_j^{-1}x_j\phi_j\in{\cal B}_T(k-1)$ and $\eta_j=\pm1$. Thus
$$
w=(\prod_{j=1}^s(\phi_j^{-1}x_j\phi_j))^{\eta_j})^{-1}\cdot\prod_{j=1}^s(y_k^{-\epsilon_k}\phi_j^{-1}x_j\phi_j y_k^{\epsilon_k}))^{\eta_j}\in F({\cal B}_T(k))
$$
The induction is finished.\\
\par
\noindent{\em Step 2.} ${\tilde i_k}: F({\cal A}_T(k))\rightarrow F({\cal B}_T(k))$ is an epimorphism,
where the homomorphism ${\tilde i_k}$ is induced by the inclusion $i_k:{\cal A}_T(k)\rightarrow F({\cal B}_T(k))$\\
\par
The proof of this step is given induction on $k$ starting with $F({\cal A}_T(0))=F({\cal B}_T(0))=F(S-T)$. Suppose that 
$F({\cal A}_T(k-1))\rightarrow F({\cal B}_T(k-1))$ is an epimorphism and consider 
${\tilde i_k}:F({\cal A}_T(k))\rightarrow F({\cal B}_T(k))$.
Let $\phi^{-1}x\phi\in{\cal B}_T(k)$, where $\phi=y_1^{\epsilon_1}\cdots y_t^{\epsilon_t}$ is a reduced
word with $t\leq k$. If $t\leq k-1$, then $\phi^{-1}x\phi\in Im\varphi_k$ by induction. Let
$\phi=y_1^{\epsilon_1}\cdots y_k^{\epsilon_k}$ be a reduced word and let $z$ denote the word
$(y_1^{\epsilon_1}\cdots y_{k-1}^{\epsilon_{k-1}})^{-1}xy_1^{\epsilon_1}\cdots y_{k-1}^{\epsilon_{k-1}}$. 
Then $\phi^{-1}x\phi=z\cdot[z,y_k^{\epsilon_k}]$. Notice that $z\in Im({\tilde i_k})$ by induction. 
 It suffices to show that 
$[w,y^{\epsilon}]\in F({\cal A}_T(k))$ for $w\in{\cal A}_T(k-1)$ for all $w\in{\cal A}_T(k-1)$, $y\in T$ and $\epsilon=\pm1$ by
the Witt-Hall identity that
$$
[ab,c]=[a,c]\cdot[[a,c],b]\cdot[b,c].
$$
We show this by second induction starting with 
$$
{\cal A}_T(1)=\{[x,y^{\epsilon}]|y\in T, x\in S-T,\epsilon=\pm1\}.
$$
 Let $w=[[x,y_1^{\epsilon_1}],\cdots],y_t^{\epsilon_t}]$ be in ${\cal A}_T(k-1)$ with $k>1$, where
$y^{\epsilon_1}_1\cdots y_t^{\epsilon_t}$ is a reduced word. Let $y\in T$ and let $\epsilon=\pm1$.
If $y_1^{\epsilon_1}\cdots y_t^{\epsilon_t}y^{\epsilon}$ is a reduced word, then 
$[w,y^{\epsilon}]\in F({\cal A}_T(k))$ by definition. Suppose that
$y_1^{\epsilon_1}\cdots y_t^{\epsilon_t}y^{\epsilon}$ is not a reduced word. Then $t>0$, $y=y_t$ and $\epsilon=-\epsilon_t$.
Let $w'$ denote $[[x,y_1^{\epsilon_1}],\cdots,],y_{t-1}^{\epsilon_{t-1}}]\in{\cal A}_T(k-2)$.
Then $w=[w',y_t^{\epsilon_t}]$. By the Witt-Hall identities, there is an equation
$$
1=[w',y_t^{-\epsilon_t}]\cdot w\cdot [w,y_t^{-\epsilon_t}].
$$
By induction, $[w',y_t^{-\epsilon_t}]\in F({\cal A}_T(k-1))\subseteq F({\cal A}_T(k))$ and so
$$[w,y_t^{-\epsilon_t}]=w^{-1}[w',y_t^{-\epsilon_t}]^{-1}\in F({\cal A}_T(k-1))\subseteq F({\cal A}_T(k)).$$
The second induction is finished and so the first induction is finished. The assertion follows.\\
\par
\noindent{\em Step 3.} ${\tilde i_k}: F({\cal A}_T(k))\rightarrow F({\cal B}_T(k))$ is an isomorphism.\\
\par
By Step 2, ${\bf Z}({\cal A}_T(k))\rightarrow{\bf Z}({\cal B}_T(k))$ is an epimorphism. Notice that
$$rank(F({\cal B}_T(k)))=|{\cal B}_T(k)|=|{\cal A}_T(k)|=rank(F({\cal A}_T(k))).$$
Thus ${\tilde i_k}:{\bf Z}({\cal A}_T(k))\rightarrow{\bf Z}({\cal B}_T(k))$ is an isomorphism and so
${\tilde i_k}: F({\cal A}_T(k))\rightarrow F({\cal B}_T(k))$ is a monomorphism. Thus
${\tilde i_k}$ is an isomorphism.\\
\par

\noindent{\em Step 4.} Since $F({\cal A}_T(k))\rightarrow F({\cal B}_T(k))$ is an isomorphism
 for each $k$,
$F({\cal A}_T)=colim_kF({\cal A}_T(k))\rightarrow F({\cal B}_T)=colim_kF({\cal B}_T(k))$ is
an isomorphism.\\
\par
Now consider the general case. By Lemma 3.6, it suffices to show that 
${\tilde i}: F({\cal A}_T)\rightarrow F({\cal B}_T)$ is an isomorphism. To check that 
$F({\cal A}_T)\rightarrow F({\cal B}_T)$ is an epimorphism. Let $w\in{\cal B}_T$, there exist finite subsets
$S'$ and $T'$ of $S$ and $T$, respectively, so that $w\in{\cal B}_{T'}$. By the special case as above,
${\tilde i}|_{F({\cal A}_{T'})}: F({\cal A}_{T'})\rightarrow F({\cal B}_{T'})$ is an isomorphism and
$w\in Im{\tilde i}|_{F({\cal A}_{T'})}$. Thus ${\tilde i}$ is an epimorphism. To check that
$F({\cal A}_T)\rightarrow F({\cal B}_T)$ is a monomorphism. Let $w\in Ker{\tilde i}$, there exist finite
subsets $S'$ and $T'$ of $S$ and $T$, respectively, so that $w\in F({\cal A}_{T'})$. Notice that
${\tilde i}|_{F({\cal A}_{T'})}$ is an isomorphism. Thus $w=1$ and the assertion follows.\\
\par
Now let's consider the intersection of kernels of projection homomorphisms. Let $S$ be a set and let $T_j$
be a subset of $S$ for $1\leq j\leq k$. Let $\pi_j:F(S)\rightarrow F(T_j)$ be the projection homomorphism
for $1\leq j\leq k$. We construct a subset ${\cal A}(T_1,\cdots,T_k)$ of $F(S)$ by induction on $k$ as follows.
$$
{\cal A}(T_1)={\cal A}_{T_1},
$$
where ${\cal A}_T$ is defined in Definition3.1. Let 
$$
T_2^{(2)}=\{w\in{\cal A}(T_1)|w=[[x,y_1^{\epsilon_1},\cdots,y_t^{\epsilon_t}]\, with\,x,y_j\in T_2\, for\,all\,j\}
$$
and define 
$$
{\cal A}(T_1,T_2)={\cal A}(T_1)_{T_2^{(2)}}.
$$
 
Suppose that ${\cal A}(T_1,T_2,\cdots,T_{k-1})$ is well defined so that all of the elements in 
${\cal A}(T_1,T_2,\cdots,T_{k-1})$ are written down as certain commutators in $F(S)$ in terms of
elements in $S$. Let $T_k^{(k)}$
be the subset of ${\cal A}(T_1,T_2,\cdots,T_{k-1})$ defined by
$$
T^{(k)}_k=\{w\in{\cal A}(T_1,T_2,\cdots,T_{k-1})|w=[x_1^{\epsilon_1},\cdots,x_l^{\epsilon_l}]\,with\,x_j\in T_k\,for\,all\,j\},
$$
where $[x_1^{\epsilon_1},\cdots,x_l^{\epsilon_l}]$ are the elements in
 ${\cal A}(T_1,T_2,\cdots,T_{k-1})$
which are written down as commutators. Then define
$$
{\cal A}(T_1,T_2,\cdots,T_k)={\cal A}(T_1,T_2,\cdots,T_{k-1})_{T_k^{(k)}}
$$

\begin{theorem}
Let $S$ be a set and let $T_j$
be a subset of $S$ for $1\leq j\leq k$. Let $\pi_j:F(S)\rightarrow F(T_j)$ be the projection homomorphism
for $1\leq j\leq k$. Then
the intersection $\bigcap_{j=1}^kKer\pi_j$ is a free group freely generated by ${\cal A}(T_1,T_2,\cdots,T_k)$.
\end{theorem}

\noindent{\em Proof.} The proof is given by induction on $k$. If $k=1$, the assertion follows from the above lemma. Suppose
that $\bigcap_{j=1}^{k-1}Ker\pi_j=F({\cal A}(T_1,T_2,\cdots,T_{k-1}))$ and consider 
$\pi_k:F(S)\rightarrow F(T_k)$. Then 
$$\bigcap_{j=1}^kKer\pi_j=Ker({\bar\pi_k}:F({\cal A}(T_1,T_2,\cdots,T_{k-1}))\rightarrow F(T_k))$$, where
${\bar\pi_k}$ is $\pi_k$ restricted to the subgroup $F({\cal A}(T_1,T_2,\cdots,T_{k-1}))$.
Let $$w=[x_1^{\epsilon_1},\cdots,x_l^{\epsilon_l}]\in{\cal A}(T_1,T_2,\cdots,T_{k-1}).$$ If $w\notin T_k^{(k)}$,
then $x_j\notin T_k$ for some $j$ and ${\bar\pi_k}(w)=1$. Thus ${\bar\pi_k}$ factors through
$F(T_k^{(k)})$, i.e, there is a homomorphism $j:F(T_k^{(k)})\rightarrow F(T_k)$ so that 
${\bar\pi_k}=j\circ\pi$, where $\pi:F({\cal A}(T_1,T_2,\cdots,T_{k-1}))\rightarrow F(T_k^{(k)})$ is the
projection homomorphism. We claim that $j:F(T_k^{(k)})\rightarrow F(T_k)$ is a monomorphism.
Consider the commutative diagram

$$
\hspace{0.0in}
\begin{array}{cccccccccc}
&F({\cal A}(T_1,T_2,\cdots,T_{k-1}))&\stackrel{{\bar\pi_k}}{\rightarrow}& F(T_k)&  
\hookrightarrow& F(S) \\  
& \uparrow &  & \uparrow j &  & \uparrow  \\  
& F(T_k^{(k)}) &\stackrel{=}{\rightarrow} & F(T_k^{(k)}) & \hookrightarrow & F({\cal A}(T_1,T_2,\cdots,T_{k-1})),  
\end{array}
$$
where $F(T_k^{(k)})\rightarrow F({\cal A}(T_1,T_2,\cdots,T_{k-1}))$ and 
$F({\cal A}(T_1,T_2,\cdots,T_{k-1}))\rightarrow F(S)$ are inclusions of subgroups. Thus 
$j:F(T_k^{(k)})\rightarrow F(T_k)$ is a monomorphism and 
$$
Ker{\bar\pi_k}=Ker(F({\cal A}(T_1,T_2,\cdots,T_{k-1}))\rightarrow F(T_k^{(k)})=F({\cal A}(T_1,T_2,\cdots,T_k)).
$$
The assertion follows.

\begin{corollary}
Let $\pi_j$ be the projection homomorphisms as in Theorem 3.5. Then the intersection subgroup
$\bigcap_{j=1}^kKer\pi_j$ equals the commutator subgroup $[[<T_1>,\cdots,<T_k>]]$ which is defined
in Notations 1.3.
\end{corollary}

\section{On the Homotopy Groups of the $3$-sphere}

In this section, we  study the Moore chain complex of the Milnor's construction $F(S^1)$. The proofs
Theorems 1.2 and 1.4 are given in this section, where Theorem 4.5 is Theorem 1.2 and Theorem 4.10 is Theorem 1.4.
 Recall that
the simplicial $1$-sphere $S^1$ is a free simplicial set generated by a $1$-simplex $\sigma$. Thus
 $S^1_0=\{*\}$, $S^1_1=\{\sigma,*\}$ and 
$S^1_{n+1}=\{*,s_n\cdots s_{i+1}s_{i-1}\cdots s_0\sigma|0\leq j\leq n\}$. Let $x_i$ denote
$s_n\cdots s_{i+1}s_{i-1}\cdots s_0\sigma$. Then $F(S^1)_{n+1}=F(x_0,x_1,\cdots, x_n)$ the free group
freely generated by $x_0,\cdots,x_n$. Let $y_i$ denote $x_{i-1}x_i^{-1}$ for $-1\leq i\leq n$, where
we put $x_{-1}=x_{n+1}=1$ in $F(S^1)_{n+1}=F(x_0,\cdots, x_n)$.
 By direct calculation, we have

\begin{lemma}
$F(S^1)_{n+1}=F(y_0,\cdots,y_n)$ with
$$
d_jy_k=\cases{y_{k-1}&$j\leq k$,\cr 1&$j=k+1$,\cr y_k&$j>k+1$,\cr}
$$
and
$$
s_jy_k=\cases{y_{k+1}&$j\leq k$,\cr y_ky_{k+1}&$j=k+1$,and\cr y_k&$j>k+1$\cr}
$$
for $0\leq j\leq n+1$, where $y_{-1}=(y_0\cdots y_{n-1})^{-1}$ in $F(S^1)_n$.
\end{lemma}

Now let $C_{n+1}$ denote the subgroup of $F(y_0,\cdots,y_n)$ generated by all of the commutators 
$[y_{i_1}^{\epsilon_1},\cdots,y_{i_t}^{\epsilon_t}]$ with $\{i_1,\cdots,i_t\}=\{0,1,\cdots,n\}$ as sets,
i.e. each $j$ $(0\leq j\leq n)$ appears in the index set $\{i_1,\cdots, i_t\}$ at least one time, where
$\epsilon_j=\pm1$ and the commutator
$[y_{i_1}^{\epsilon_1},\cdots,y_{i_t}^{\epsilon_t}]$ runs over all of the commutator bracket 
arrangements of weight $t$ for $y_{i_1}^{\epsilon_1},\cdots,y_{i_t}^{\epsilon_t}$.

\begin{lemma}
The group $C_{n+1}$ is a subgroup of $NF(S^1)_{n+1}$,i.e $C_{n+1}\subseteq \cap_{j\neq 0}Ker(d_j)$.
\end{lemma}

\noindent{\em Proof.} Notice that $d_jy_{j-1}=1$  for $1\leq j\leq n+1$. Since
 $\{y_{i_1},\cdots, y_{i_t}\}=\{y_0,y_1,\cdots,y_n\}$,
 $d_j[y_{i_1}^{\epsilon_1},\cdots,y_{i_t}^{\epsilon_t}]=1$ for each $j>0$. The assertion follows.\\
\par

\begin{theorem}
$NF(S^1)_{n+1}=C_{n+1}$.
\end{theorem}

\noindent{\em Proof.} For $1\leq j\leq n+1$, let $S=\{y_0,y_1,\cdots,y_n\}$ and let 
$T_j=\{y_0,\cdots,{\hat y_j}\cdots,y_n\}$. By Lemma 4.1, there is a commutative diagram

$$
\hspace{0.0in}
\begin{array}{cccccc}
& F(S) & \stackrel{\pi_j}{\rightarrow} & F(T_{j-1}) \\  
& \downarrow d_j &  & {\bar d_j}\downarrow \cong   \\  
& F(y_0,\cdots,y_{n-1}) &\stackrel{=}{\rightarrow} & F(y_0,\cdots,y_{n-1}),\\  
\end{array}
$$
where 
$$
{\bar d_j}(y_k)=\cases{y_k&$0\leq k\leq j-2,$\cr y_{k-1}&$j\leq k\leq n.$\cr}
$$
Thus $Kerd_j=Ker\pi_j$ and $NF(S^1)_{n+1}=\bigcap_{j=1}^{n+1}Kerd_j=\bigcap_{j=1}^{n+1}Ker\pi_j$. 
By Theorem 3.5, $NF(S^1)_{n+1}=F({\cal A}(T_0,T_1,\cdots,T_n))$, where the notation
${\cal A}(T_0,T_1,\cdots,T_n)$ is given in Section 3. To check that 
$$
F({\cal A}(T_0,T_1,\cdots,T_n))\subseteq C_{n+1},
$$
it suffices to show that 
${\cal A}(T_0,T_1,\cdots,T_n)\subseteq C_{n+1}$. This will follow from the following statement.\\
\par
\noindent{\em Statement.} For each $0\leq j\leq n$ and 
$w=[y_{i_1}^{\epsilon_1},y_{i_2}^{\epsilon_2},\cdots,y_{i_t}^{\epsilon_t}]\in{\cal A}(T_0,T_1,\cdots,T_j)$,
$\{y_0,y_1,\cdots,y_j\}\subseteq\{y_{i_1},y_{i_2},\cdots,y_{i_t}\}$.\\
\par
We show this statement by induction on $j$. Note that $F(T_0)=F(y_1,\cdots,y_n)$. For $j=0$, 
$$
{\cal A}(T_0)=\{[[y_0,y_{i_1}^{\epsilon_1}],\cdots,],y_{i_t}^{\epsilon_t}]|i_s>0,
y_{i_1}^{\epsilon_1}\cdots y_{i_t}^{\epsilon_t}\,a\,reduced\,word\, in\,F(T_0)\}.
$$
Thus the assertion holds for $j=0$. Suppose that the assertion holds for $j-1$ with $j\leq n$.
Notice that $T_j=\{y_0,\cdots,{\hat y_j},\cdots, y_n\}$. Thus
$$
T_j^{(j)}=\{w\in{\cal A}(T_0,T_1,\cdots,T_{j-1})|w=[y_{i_1}^{\epsilon_1},\cdots,y_{i_t}^{\epsilon_t}]\,
with\,y_j\notin\{y_{i_1},\cdots,y_{i_t}\}\}.
$$
and so $y_j\in\{y_{i_1},\cdots,y_{i_t}\}$ for  $w=[y_{i_1}^{\epsilon_1},\cdots,y_{i_t}^{\epsilon_t}]\in{\cal A}(T_0,T_1,\cdots,T_{j-1})-T_j^{(j)}$.
 Hence, by induction, 
$$
\{y_0,y_1,\cdots,y_j\}\subseteq\{y_{i_1},y_{i_2},\cdots,y_{i_t}\}
$$
for
$w=[y_{i_1}^{\epsilon_1},\cdots,y_{i_t}^{\epsilon_t}]\in{\cal A}(T_0,T_1,\cdots,T_{j-1})-T_j^{(j)}$.
Notice that 
$${\cal A}(T_0,T_1,\cdots,T_j)={\cal A}(T_0,T_1,\cdots,T_{j-1})_{T_j^{(j)}}.$$ 
Thus
$$
\{y_0,y_1,\cdots,y_j\}\subseteq\{y_{i_1},y_{i_2},\cdots,y_{i_t}\}
$$
for each
$w=[y_{i_1}^{\epsilon_1},\cdots,y_{i_t}^{\epsilon_t}]\in{\cal A}(T_0,T_1,\cdots,T_{j})$. The induction
is finished and the theorem follows.

\begin{corollary}
 ${\cal A}(T_0,T_1,\cdots,T_n)$ is a set of free generators for $NF(S^1)_{n+1}$.
\end{corollary}

\begin{theorem}[Theorem 1.2]
For $n\geq1$, $\pi_{n+2}(S^3)$ is  isomorphic to the center of the group with generators $y_0,y_1,\cdots,y_n$
and relations $[y_{i_1}^{\epsilon_1},y_{i_2}^{\epsilon_2},\cdots,y_{i_t}^{\epsilon_t}]$ with
$$
\{i_1,i_2,\cdots,i_t\}=\{-1,0,1,\cdots,n\}
$$
 as sets, where
$\epsilon_j=\pm1$,$y_{-1}=(y_0y_1\cdots y_n)^{-1}$ and the commutator bracket $[\cdots]$ runs over all
bracket arrangements of weight $t$ for each $t$.
\end{theorem}

\noindent{\em Proof.} Notice that $\pi_{n+2}(S^3)\cong\pi_{n+1}F(S^1)$ for $n\geq1$ and 
${\cal B}_{n+1}=d_0(NF(S^1)_{n+2})$. By Lemma 4.1 and Theorem 4.3, ${\cal B}_{n+1}$ is generated by 
$[y_{i_1}^{\epsilon_1},y_{i_2}^{\epsilon_2},\cdots,y_{i_t}^{\epsilon_t}]$ so that 
$\{y_{i_1},y_{i_2},\cdots,y_{i_t}\}=\{y_0,y_1,\cdots,y_n\}$, where $\epsilon_j=\pm1$. By Proposition 2.8,
it suffices to check that $\pi_2(F(S^1))\cong Z(F(S^1)_2/{\cal B}_2)=Z(F(y_0,y_1)/{\cal B}_2)$. By 
Example 2.13, ${\cal Z}_2=\Gamma^2F(S^1)_2=\Gamma^2F(y_0,y_1)$ and 
$\Gamma^3F(y_0,y_1)\subseteq{\cal B}_2$. Now , by the construction of ${\cal B}_2$, 
${\cal B}_2\subseteq\Gamma^3F(y_0,y_1)$ and therefore 
$$\pi_2(F(S^1))\cong{\cal Z}_2/{\cal B}_2=\Gamma^2F(y_0,y_1)/\Gamma^3F(y_0,y_1)=Z(F(y_0,y_1)/{\cal B}_2).$$
The assertion follows.

\begin{remark}
 The relations in above theorem are not minimal, i.e, many of them can be cancelled out.
\end{remark}

\begin{proposition}
Let ${C'}_{n+1}$ be the subgroup of $F(S^1)_{n+1}$ generated by all commutators given by
$[\cdots[y_{i_1}^{\epsilon_1},y_{i_2}^{\epsilon_2}],\cdots,],y_{i_t}^{\epsilon_t}]$ with
$\{i_1,\cdots,i_t\}=\{0,1,\cdots,n\}$ and $\epsilon_j=\pm1$. Then 
$$
NF(S^1)_{n+1}/\Gamma^s\cap NF(S^1)_{n+1}\cong C'_{n+1}/\Gamma^s\cap C'_{n+1}
$$
for each $s$, where $\Gamma^s=\Gamma^sF(S^1)_{n+1}$ is the $s$-term in the lower central series of $F(S^1)_{n+1}$.
\end{proposition}

\noindent{\em Proof.} Notice that $C'_{n+1}\subseteq NF(S^1)_{n+1}$. The induced homomorphism
$f_s:C'_{n+1}/\Gamma^s\cap C'_{n+1}\rightarrow NF(S^1)_{n+1}/\Gamma^s\cap NF(S^1)_{n+1}$ is a monomorphism.
We check that $f$ is an epimorphism. It suffices to show that,
for each $w\in NF(S^1)$, there exists a sequence of elements
$\{x_j\}$ so that $x_j\in\Gamma^j\cap C'_{n+1}$ and $wx_1x_2\cdots x_s\in\Gamma^{s+1}$ for each $s$.
In fact, if this statement holds, 
$wx_1x_2\cdots x_{s-1}\equiv 1$ $mod$ $\Gamma^s\cap NF(S^1)_{n+1}$ 
for each $s$ and $w=(wx_1\cdots x_{s-1})\cdot(x_1x_2\cdots x_{j-1})^{-1}\in C'_{n+1}$
$mod\,\Gamma^s$.
 Now we 
construct $x_j$ by induction, which depends on $w$. Notice that, for $n\geq1$, 
$NF(S^1)_{n+1}\in\Gamma^{n+1}\subseteq\Gamma^2$. Choose $x_j=1$ for $j\leq n$.
 Suppose that there are
$x_1,\cdots,x_{s-1}$ so that $x_j\in\Gamma^j\cap C'_{n+1}$ and $wx_1\cdots x_{s-1}\in\Gamma^s$.
Since $C'_{n+1}\subseteq NF(S^1)_{n+1}$, $wx_1\cdots x_{s-1}\in\Gamma^s\cap NF(S^1)_{n+1}$ and
$d_j(wx_1\cdots x_{s-1})=1$ for $j>1$. Let
$\pi:\Gamma^s\rightarrow\Gamma^s/\Gamma^{s+1}$ be the quotient homomorphism.
$\pi(wx_1\cdots x_{s-1})$ is a linear combination of basic Lie products. We claim that
$\pi(wx_1\cdots x_{s-1})$ is a linear combination of basic Lie products in  which each $y_j$ appears
in the Lie product for $0\leq j\leq n$. If not, there exists $j$ so that 
$\pi(wx_1\cdots x_{s-1})=b+c$, where $b$ is a nontrivial linear combination of basic Lie products in which
$y_j$ does not appear and $c$ is a linear combination of basic Lie products in which $y_j$ appears.
Now the face homomorphism $d_{j+1}:F(y_0,\cdots,y_n)\rightarrow F(y_0,\cdots,y_{n-1})$
 induces a homomorphism
$$
{\bar d}_{j+1}:\Gamma^sF(y_0,\cdots,y_n)/\Gamma^{s+1}F(y_0,\cdots,y_n)\rightarrow\Gamma^sF(y_0,\cdots,y_{n-1})/\Gamma^{s+1}F(y_0,\cdots,y_{n-1})
$$
and
 $$1={\bar d}_{j+1}\pi(wx_1\cdots x_{s-1})={\bar d}_{j+1}(b)+{\bar d}_{j+1}(c)={\bar d}_{j+1}(b).$$
Notice that 
$$d_{j+1}|_{F(y_0,\cdots,y_{j-1},y_{j+1},\cdots,y_n)}:F(y_0,\cdots,y_{j-1},y_{j+1},\cdots,y_n)\rightarrow F(y_0,\cdots,y_{n-1})$$
is an isomorphism. Thus $b=1$. This contradicts to that $b$ is a nontrivial linear combination of a basis.
Thus $\pi(wx_1\cdots x_{s-1})$ is a linear combination of basic Lie products in which all of the $y_j$
appear. By Theorem 5.12 in [MKS,pp.337], there exists $x_s$ in $C'_{n+1}$ so that $\pi(wx_1\cdots x_s)=1$, or
$wx_1\cdots x_s\in\Gamma^{s+1}$. The induction is finished now. \\
\par
By Corollary 3.6, we have

\begin{theorem}
In the free group $F(y_0,\cdots,y_n)$, $NF(S^1)_{n+1}=[[<y_0>,\cdots,<y_n>]]$ and 
${\cal B}F(S^1)_{n+1}=[[<y_{-1}>,<y_0>,\cdots,<y_n>]]$.
\end{theorem}
Thus Theorem 4.5 can be rewritten as follows.
\begin{theorem} In the free group $F(y_0,\cdots,y_n)$ for $n\geq1$, the center
$$
Z(F(y_0,\cdots,y_n)/[[<y_{-1}>,<y_0>,\cdots,<y_n>]])\cong\pi_{n+2}(S^3)
$$   
\end{theorem}
\par
By Lemma 4.1 and Proposition 3.5, $Kerd_0=<x_0>=<y_{-1}>$ and therefore 
$$
{\cal Z}_{n+1}=[[<y_0>,\cdots,<y_n>]]\cap<y_{-1}>.
$$ Thus we have
\begin{theorem}[Theorem 1.4]
In the free group $F(y_0,\cdots,y_n)$ with $n\geq1$,
$$
[[<y_0>,\cdots,<y_n>]]\cap<y_{-1}>]]/[[<y_{-1}>,<y_0>,\cdots,<y_n>]]\cong\pi_{n+2}(S^3).
$$
\end{theorem}

\section{On the Homotopy Groups of $\Sigma K(\pi,1)$}
In this section, we give group theoretical descriptions for $\pi_*(\Sigma K(\pi,1))$ for any 
group $\pi$. The proof of Theorem 1.5 is given in this section, where Theorem 5.9 is Theorem 1.5. We will use the notations defined in Section 3.
 First we extend our description for $\pi_*(S^2)$ to the case $\pi_*(\vee_{\alpha\in J}S^2)$.
Recall that $(\vee_{\alpha\in J}S^1)_0=*$, $(\vee_{\alpha\in J}S^1)_1=\{\sigma_{\alpha},*|\alpha\in J\}$
and 
$(\vee_{\alpha\in J}S^1)_{n+1}=\{s_n\cdots{\hat s_i}\cdots s_0\sigma_{\alpha},*|\alpha\in J,0\leq i\leq n\}$.
Let $x^{(\alpha)}_i$ denote $s_n\cdots{\hat s_i}\cdots s_0\sigma_{\alpha}$. Then
$$
F(\vee_{\alpha\in J}S^1)_{n+1}=F(x^{(\alpha)}_0,x^{(\alpha)}_1,\cdots,x^{(\alpha)}_n|\alpha\in J).
$$
Let $y^{(\alpha)}_j=x^{(\alpha)}_j\cdot x^{(\alpha)}_{j+1}$ for $0\leq j\leq n-1$ and $y^{(\alpha)}_n=x^{(\alpha)}_n$.
By Lemma 4.1, we have
\begin{lemma}
$F(\vee_{\alpha\in J}S^1)_{n+1}=F(y^{(\alpha)}_j|0\leq j\leq n,\alpha\in J)$ with
$$
d_j(y^{(\alpha)}_k)=\cases{y^{(\alpha)}_{k-1}& $j\leq k$, \cr 1& $j=k+1$, \cr
y^{(\alpha)}_k & $j>k+1$, \cr}  
$$
and 
$$
s_j(y^{(\alpha)}_k)=\cases{y^{(\alpha)}_{k+1}& $j\leq k$, \cr y^{(\alpha)}_k\cdot y^{(\alpha)}_{k+1}& $j=k+1$, \cr
y^{(\alpha)}_k & $j>k+1$, \cr}  
$$
for $0\leq j\leq n+1$, where $y^{(\alpha)}_{-1}=(y^{(\alpha)}_0\cdots y^{(\alpha)}_{n-1})^{-1}$ in
$F(\vee_{\alpha\in J}S^1)$
\end{lemma}

Let $C^J_{n+1}$ denote the subgroup of $F(\vee_{\alpha\in J}S^1)_{n+1}$ generated by all of the commutators 
$[y^{(\alpha_1)\epsilon_1}_{i_1},\cdots, y^{(\alpha_t)\epsilon_t}_{i_t}]$ with
 $\{i_1,\cdots,i_t\}=\{0,1,\cdots,n\}$ as sets, where $\epsilon_j=\pm1$, $\alpha_j\in J $ and the 
commutator bracket $[\cdots]$ runs over all of the commutator bracket arrangements of weight $t$ for
each $t$.

\begin{lemma}
$C^J_{n+1}\subseteq NF(\vee_{\alpha\in J}S^1)_{n+1}$.
\end{lemma}

\noindent{\em Proof:} For each $1\leq j\leq n+1$, there exists some $i_s=j-1$. Thus
 $d_j(y^{(\alpha_s)\epsilon_s}_{i_s})=1$ for some $i_s$
and therefore
 $$
d_j([y^{(\alpha_1)\epsilon_1}_{i_1},\cdots, y^{(\alpha_t)\epsilon_t}_{i_t}])=1
$$
for each $j>0$. 
 The assertion follows.

\begin{lemma}
For each $1\leq j\leq n+1$, 
$$
Kerd_j\cap F(\vee_{\alpha\in J}S^1)_{n+1}=<y^{(\alpha)}_{j-1}|\alpha\in J>,
$$
the normal subgroup generated by $y_{j-1}^{(\alpha)}$ with $\alpha\in J$.
\end{lemma}
\noindent{\em Proof:} By the definition of $d_j$, there is a commutative diagram
$$
\hspace{0.0in}
\begin{array}{cccccc}
& F(y^{(\alpha)}_j|0\leq j\leq n,\alpha\in J) & \stackrel{p}{\rightarrow} & F(y^{(\alpha)}_0\cdots{\hat y^{(\alpha)}_{j-1}}\cdots y_n^{(\alpha)}|\alpha\in J) \\  
& d_j\downarrow &  &\cong \downarrow {\bar d_j} \\  
& F(y^{(\alpha)}_j|0\leq j\leq n-1,\alpha\in J)&\stackrel{=}{\rightarrow} & F(y^{(\alpha)}_j|0\leq j\leq n-1,\alpha\in J) \\  
\end{array}
$$
where $p$ is the projection and
$$
{\bar d_j}y^{(\alpha)}_k=\cases{y^{(\alpha)}_{k-1}& $j\leq k$,\cr y^{(\alpha)}_k& $j>k+1$.\cr}
$$
The assertion follows.
 
\begin{theorem}
Let $C^J_{n+1}$ be defined as above. Then
$$
NF(\vee_{\alpha\in J}S^1)_{n+1}=C^J_{n+1}
$$
\end{theorem}
\noindent{\em Proof:}
 By lemma 5.1, each $d_j$ with $j>0$ is a projection homomorphism. Thus, by Theorem 3.5, 
$$
NF(\vee_{\alpha\in J}S^1)_{n+1}=\bigcap_{j=1}^{n+1}Kerd_j=F({\cal A}(T_0,T_1,\cdots,T_n))
$$
where
$T_j=\{y^{(\alpha)}_0,\cdots,{\hat y^{(\alpha)}_j},\cdots,y^{(\alpha)}_n|\alpha\in J\}$. It suffices to show that
$${\cal A}(T_0,T_1,\cdots,T_n)\subseteq C^J_{n+1}.$$ This follows from the next lemma.
\begin{lemma}
 For $0\leq j\leq n$, let $W=[y^{(\alpha_1)\epsilon_1}_{i_1},\cdots,y^{(\alpha_t)\epsilon_t}_{i_t}]\in {\cal A}(T_0,T_1,\cdots,T_j)$.
Then
$$
\{0,1,\cdots,j\}\subseteq\{i_1,\cdots,i_t\}.
$$
\end{lemma}
\noindent{\em Proof:} The proof is given by induction on $j$  for $0\leq j\leq n$. Notice that, by construction, each 
element in ${\cal A}(T_0,T_1,\cdots,T_j)$ is  written as a certain commutator. If $j=0$, then
$$
{\cal A}(T_0)=\{y^{(\alpha)}_0,[[y^{(\alpha)}_0,y^{(\alpha_1)\epsilon_1}_{i_1}],\cdots,],y^{(\alpha_t)\epsilon_t}_{i_t}]\}
$$
where $\alpha,\alpha_j\in J$, $\epsilon_j=\pm1$ and $y^{(\alpha_1)\epsilon_1}_{i_1}\cdots y^{(\alpha_t)\epsilon_t}_{i_t}$
runs over all of the reduced words $\not=1$ in $F(y^{(\alpha)}_j|\alpha\in J,1\leq j\leq n)$.
Thus the assertion holds for $j=0$. Suppose that the the assertion holds for $j-1$ with $j\leq n$.
Recall that 
$$
T^{(j)}_j=\{W\in{\cal A}(T_0,T_1,\cdots,T_{j-1}|W=[y^{(\alpha_1)\epsilon_1}_{i_1},y^{(\alpha_2)\epsilon_2}_{i_2},\cdots,y^{(\alpha_t)\epsilon_t}_{i_t}],\,with\,y^{(\alpha_s)}_{i_s}\in T_j\}.
$$
Notice that $y^{(\alpha_s)}_{i_s}\in T_j\Longleftrightarrow i_s\not=j$. Thus,
for $W=[y^{(\alpha_1)\epsilon_1}_{i_1},y^{(\alpha_2)\epsilon_2}_{i_2},\cdots,y^{(\alpha_t)\epsilon_t}_{i_t}]\in{\cal A}(T_0,\cdots,T_{j-1})-T^{(j)}_j$,
$j\in\{i_1,\cdots,i_t\}$. By induction, $\{0,\cdots,j-1\}\subseteq\{i_1,\cdots,i_t\}$. Hence
$$\{0,1,\cdots,j\}\subseteq\{i_1,i_2,\cdots,i_t\}$$
for any 
$W=[y^{(\alpha_1)\epsilon_1}_{i_1},y^{(\alpha_2)\epsilon_2}_{i_2},\cdots,y^{(\alpha_t)\epsilon_t}_{i_t}]\in{\cal A}(T_0,\cdots,T_{j-1}-T^{(j)}_j$.
Recall that, by construction, 
$$
{\cal A}(T_0,\cdots, T_j)={\cal A}_{T^{(j)}_j}. 
$$
$$\{0,1,\cdots,j\}\subseteq\{i_1,i_2,\cdots,i_t\}$$
for any 
$W=[y^{(\alpha_1)\epsilon_1}_{i_1},y^{(\alpha_2)\epsilon_2}_{i_2},\cdots,y^{(\alpha_t)\epsilon_t}_{i_t}]\in{\cal A}(T_0,\cdots,T_{j})$.
This completes the proof.

\begin{theorem}
$\pi_{n+2}(\vee_{\alpha\in J}S^2)$ is isomorphic to the center of the group with generators
$$
y^{(\alpha)}_0,y^{(\alpha)}_1,\cdots,y^{(\alpha)}_n
$$
for $\alpha\in J$. and relations
$$
[y^{(\alpha_1)\epsilon_1}_{i_1}, y^{(\alpha_2)\epsilon_2}_{i_2},\cdots,y^{(\alpha_t)\epsilon_t}_{i_t}]
$$ 
with $\{i_1,i_2,\cdots,i_t\}=\{-1,0,1,\cdots,n\}$ as sets, where the indices $i_j$ can be repeated, 
$\epsilon_j=\pm1$, $\alpha_j\in J$ and the commutator bracket runs over all of the commutator bracket
arrangements of weight $t$ for each $t$.
\end{theorem}

\noindent{\em Proof:} Notice that $\pi_{n+2}(\vee_{\alpha\in J} S^2)\cong\pi_{n+1}(F(\vee_{\alpha\in J}S^1))$.
 By the above theorem, ${\cal B}_{n+1}$ is generated by
$$
[y^{(\alpha_1)\epsilon_1}_{i_1}, y^{(\alpha_2)\epsilon_2}_{i_2},\cdots,y^{(\alpha_t)\epsilon_t}_{i_t}]
$$ 
By Theorem 2.8, the assertion holds for $n\geq1$. For $n=0$, ${\cal B}_1=\Gamma^2(F(y^{(\alpha)}_0|\alpha\in J))$
and 
$$
F(y^{(\alpha)}_0|\alpha\in J)/{\cal B}_1\cong \oplus_{\alpha\in J}{\bf Z}\pi_2(\vee_{\alpha\in J}S^2).
$$ 
The assertion holds for all of $n$.\\
\par
For the general case, we need a simplicial group construction.
\begin{definition}
Let $G$ be a simplicial group and let $X$ be a pointed simplicial set with a point $*$. The simplicial group
{\bf $F^G(X)$} is defined by setting
$$
F^G(X)_n=\coprod_{x\in X_n}(G_n)_x,
$$
the free product, modulo the relations $(G_n)_*$, where $(G_n)_x$ is a copy of $G_n$. The faces and
degeneracies homomorphisms in $F^G(X)$ is given in the canonical way by the universal property of the
coproduct in the category of groups and group homomorphisms.
\end{definition}

\begin{lemma}[Ca, Theorem 9, pp.88] 
Let $G$ be a simplicial group and let $X$ be a pointed simplicial set. Then the goemetric realization
$|F^G(X)|$ is homotopy equivalent to $\Omega(B|G|\wedge|X|)$.
\end{lemma}

A generalization of this lemma by using fibrewise simplicial groups is given in [Wu1].

\begin{theorem}[Theorem 1.5]
Let $\pi$ be any group and let $\{x^{(\alpha)}|\alpha\in J\}$ be a set of generators for $\pi$. Then, for
$n\not=1$, $\pi_{n+2}(\Sigma K(\pi,1))$ is  isomorphic to the center of the quotient group
of the free product 
$${\coprod^{groups}}_{0\leq j\leq n}(\pi)_j$$
modulo the relations
$$
[y^{(\alpha_1)\epsilon_1}_{i_1},y^{(\alpha_2)\epsilon_2}_{i_2},\cdots,y^{(\alpha_t)\epsilon_t}_{i_t}]
$$
with $\{i_1,i_2,\cdots,i_t\}=\{-1,0,1,\cdots,n\}$ as sets,
where $(\pi)_j$ is a copy of $\pi$ with generators $\{x^{(\alpha)}_j|\alpha\in J\}$, $\epsilon_j=\pm1$,
$y^{(\alpha)}_{-1}={x^{(\alpha)}_0}^{-1}$, $y_j=x^{(\alpha)}_j{x^{(\alpha)}_{j+1}}^{-1}$ for $1\leq j\leq n-1$ and
 $y^{(\alpha)}_n=x^{(\alpha)}_n$, the commutator bracket
$[\cdots]$ runs over all of the commutator bracket arrangment of weight $t$ for each $t$.
\end{theorem}
\noindent{\em Proof:} Since $\{x^{(\alpha)}|\alpha\in J\}$ is a set of generators for $\pi$,
$F(x^{(\alpha)}|\alpha\in J)\rightarrow\pi $ is an epimorphism. Thus
$$
F(\vee_{\alpha\in J}S^1)\cong F^{F(x^{(\alpha)}|\alpha\in J)}(S^1)\rightarrow F^{\pi}(S^1)
$$ is an epimorphism.
Hence 
$$
NF(\vee_{\alpha\in J}S^1)\rightarrow NF^{\pi}(S^1)
$$
and
$$
{\cal B}F(\vee_{\alpha\in J}S^1)\rightarrow {\cal B}F^{\pi}(S^1)
$$
are epimorphisms.
The assertion follows from Theorem 5.6, Lemma 5.7 and Proposition 2.8.

\section{Applications}
In this section, we consider Cohen's $K$-construction. Our descriptions for the homotopy groups of
$3$-sphere gives a calculation for the $K$-construction of $1$-sphere.

\begin{definition}
Let $X$ be set. The group ${\bf K(X)}$ is defined to the the quotient group of the free group $F(X)$
modulo the normal subgroup generated by all of the 
commutators $[[x_1,x_2],\cdots,],x_t]$ with $x_i\in X$ and $x_i=x_j$ for some $1\leq i<j\leq t$ 
Now let $S$ be a pointed simplicial set. The simplicial group ${\bf K(S)}$ is defined to be the quotient
simplicial group of $F(S)$ modulo the normal simplicial subgroup generated by all of the 
commutators $[[x_1,x_2],\cdots,],x_t]$ with $x_i\in S$ and $x_i=x_j$ for some $1\leq i<j\leq t$,
where $F(S)$ is the Milnor's $F(K)$-construction for the simplicial set $S$.
\end{definition}
\begin{definition}
The group $Lie(n)$ is the 
elments of weight $n$ in the Lie algebra $Lie(x_1,x_2,\cdots,x_n)$ which is the quotient Lie algebra of
the free Lie algebra $$L(x_1,x_2,\cdots,x_n)$$ over ${\bf Z}$ modulo the two sided Lie ideal generated by the Lie
elements $[[x_{i_1},x_{i_2}],\cdots,],x_{i_t}]$ with $i_l=i_k$ for some $1\leq l<k\leq t$.
\end{definition} 

The following lemmas are due to Fred Cohen.
\begin{lemma}[Co]
$\Gamma^qK(x_1,x_2,\cdots,x_n)=\{1\}$ for $q\geq n$ and 
$$\Gamma^nK(x_1,x_2,\cdots,x_n)\cong Lie(n),$$
where $\Gamma^qG$ is the $q$-th term in the lower central series for the group $G$  starting with $\Gamma^1G=G$.
\end{lemma}
\begin{lemma}[Co]
In the group $K(x_1,x_2,\cdots, x_n)$, the normal subgroup grnerated by $x_j$ is abelian for each $1\leq j\leq n$.
\end{lemma}
\begin{lemma}[Co]
The set  $\{[x_1,x_{\sigma(2)},\cdots,x_{\sigma(n)}]|\sigma\in\Sigma_{n-1}\}$
 is a ${\bf Z}$-basis for $Lie(n)$, where $\Sigma_{n-1}$ acts on $\{2, 3,\cdots,n\}$ by permutation.
\end{lemma}

Recall that
the simplicial $1$-sphere $S^1$ is a free simplicial set generated by a $1$-simplex $\sigma$. More
precisely, $S^1_0=\{*\}$, $S^1_1=\{\sigma,*\}$ and 
$S^1_{n+1}=\{*,s_n\cdots{\hat s_j} \cdots s_0\sigma|0\leq j\leq n\}$. Let
$x_i$ denote $s_n\cdots {\hat s_j}\cdots s_0\sigma$. Then 
\begin{lemma} 
The face functions $d_i:S^1_{n+1}\rightarrow S^1_n$ and the degenarate functions 
$s_i:S^1_{n+1}\rightarrow S^1_{n+2}$ are as follows:
$$
d_ix_j=\cases{x_j&$j<i$\cr x_{j-1}&$j\geq i$\cr}
$$
and
$$
s_ix_j=\cases{x_j&$j<i$\cr x_{j+1}&$j\geq i$\cr},
$$
where we put $x_{-1}=*$ and $x_n=*$ in $S^1_n$.
\end{lemma}

\begin{theorem} 
$\pi_n(K(S^1))$ is isomorphic to $Lie(n)$
\end{theorem}
\noindent{\em Proof:} Let $\pi:F(S^1)\rightarrow K(S^1)$ be the quotient homomorphism. Then
$NF(S^1)\rightarrow NK(S^1)$ is an epimorphism. Recall that $NF(S^1)_{n+1}$ is generated by all of
the commutators 
$$
[y_{i_1},y_{i_2},\cdots,y_{i_t}]
$$
so that 
$$\{i_1,i_2,\cdots,i_t\}=\{0,1,\cdots,n\}$$ by Theorem 4.3. Thus
$
NF(S^1)_{n+1}\subseteq \Gamma^{n+1}F(S^1)_{n+1}
$
and therefore
$$
NK(S^1)_{n+1}\subseteq\Gamma^{n+1}K(S^1)_{n+1}.
$$
Notice that $K(S^1)_{n+1}\cong K(x_0,x_1,\cdots,x_n)$. Thus $\Gamma^{n+1}K(S^1)_{n+1}\cong Lie(n+1)$ and
$\Gamma^{n+1}K(S^1)_n=\{1\}$. Thus 
$$
d_j|_{\Gamma^{n+1}K(S^1)_{n+1}}:\Gamma^{n+1}K(S^1)_{n+1}\rightarrow K(S^1)_n
$$
is trivial for each $j\geq0$. And therefore
$$
NKS^1)_{n+1}=\Gamma^{n+1}K(S^1)_{n+1}\cong Lie(n+1)
$$
with
$
d_0:NK(S^1)_{n+1}\rightarrow NK(S^1)_n
$ is trivial. the assertion follows.

\begin{corollary}
Let $\pi:F(S^1)\rightarrow K(S^1)$ be the quotient simplicial homomorphism. Then
$$
\pi_*:\pi_n(F(S^1))\rightarrow \pi_n(K(S^1))
$$
is an isomorphism for $n=1,2$ and zero for $n>2$.
\end{corollary}
Now we consider the Samelson product in $\pi_*(K(S^1))$. Let $x_j$ denote $s_n\cdots{\hat s_j}\cdots s_0\sigma$ in
$S^1_{n+1}$. The following lemma follows directly from Lemma 6.6.

\begin{lemma}
Let $I=(i_1,i_2,\cdots,i_m)$ be a sequence with $i_1<i_2<\cdots<i_m$. Then
$ s_I:S^1_{n+1}-*\rightarrow S^1_{n+m+1}-*$ is the composite
$$
\{x_0,x_1,\cdots,x_n\}\stackrel{{\bar s_I}}{\rightarrow}\{x_0,x_1,\cdots,{\hat x_{i_1}},\cdots,{\hat x_{i_2}},\cdots,{\hat x_{i_m}},\cdots,x_{n+m}\}
$$
$$
\hookrightarrow\{x_0,x_1,\cdots,x_{n+m}\}. 
$$
where $s_I=s_{i_m}\cdots s_{i_1}$ and ${\bar s_I}$ is the order preserving isomorphism.
\end{lemma}
Recall that, for 
$\sigma\in\pi_n(G)$
 and $\tau\in\pi_m(G)$, the Samelson product [C1] is defined to be
$$
<\sigma,\tau>=\prod_{(a,b)}[s_b\sigma, s_a\tau]^{sign(a,b)},
$$
where $G$ is a simplicial group, $(a,b)=(a_1,\cdots,a_n,b_1,\cdots,b_m)$ runs over all shuffles of 
$(0,1,\cdots,m+n-1)$,i.e. all permutations, so that $a_1<a_2<\cdots<a_n$, $b_1<b_2<\cdots<b_m$, $sign(a,b)$
is the sign of the permutation $(a,b)$, the order of the product $\prod$ is right lexicographic on $a$ and
$s_a=s_{a_n}\cdots s_{a_1}$.

\begin{proposition}
The Samelson product in $\pi_*(K(S^1))$ is as follows:
$$
<[x_{\sigma(0)},x_{\sigma(1)},\cdots,x_{\sigma(n)}],[x_{\tau(0)},x_{\tau(1)},\cdots,x_{\tau(m)}]>
$$
$$
=\sum_{(I,J)}sign(I,J)[[x_{i_{\sigma(0)}},\cdots,x_{i_{\sigma(n)}}],[x_{j_{\tau(0)}},\cdots,x_{j_{\tau(m)}}]]
$$
for the commutators 
$$
[x_{\sigma(0)},x_{\sigma(1)},\cdots,x_{\sigma(n)}]\in\pi_{n+1}(KS^1))\cong Lie(n+1)
$$
and
$$
[x_{\tau(0)},x_{\tau(1)},\cdots,x_{\tau(m)}]\in\pi_{m+1}(K(S^1))\cong Lie(m+1)
$$
where 
$$(I,J)=(i_0,i_1,\cdots,i_n,j_0,j_1,\cdots,j_m)$$
 runs over all shuffles of 
$(0,1,\cdots,m+n+1)$  so that $i_0<i_1<\cdots<i_n$, $j_0<j_1<\cdots<j_m$, $sign(I,J)$
is the sign of the permutation $(I,J)$, $\sigma\in\Sigma_{n+1}$ acts on $\{0,1,\cdots,n\}$ and
$\tau\in\Sigma_{m+1}$ acts on $\{0,1,\cdots,m\}$
\end{proposition}

\noindent{\em Proof:} Notice that 
$
\{x_0,\cdots,{\hat x_{j_0}},\cdots,{\hat x_{j_m}},\cdots,x_{n+m+1}\}=\{x_{i_0},\cdots,x_{i_n}\}
$
and
$
{\bar s_J}:\{x_0,\cdots,x_n\}\rightarrow\{x_{i_0},\cdots,x_{i_n}\}
$
is an ordered isomorphism.
$$
s_J([x_{\sigma(0)},x_{\sigma(1)},\cdots,x_{\sigma(n)}])=[x_{i_{\sigma(0)}},\cdots,x_{i_{\sigma(n)}}]
$$
and
$$
s_I([x_{\tau(0)},x_{\tau(1)},\cdots,x_{\tau(m)}])=[x_{j_{\tau(0)}},\cdots,x_{j_{\tau(m)}}].
$$
The assertion follows.

\begin{definition}
A simplicial group is {\bf minimal} if it is also a minimal simplicial set.
\end{definition}
Recall that a simplicial group $G$ is minimal if and only if the Moore chain complex $NG$ is minimal [C2].
\begin{theorem}
The simplicial group $K(S^1)$ is the universal minimal simplicial quotient simplicial group of $F(S^1)$ in
the following sense:\\
(1).  $K(S^1)$ is a minimal simplicial group.\\
(2). Let $G$ be a minimal simplicial group. Then every simplicial homomorphism $f:F(S^1)\rightarrow G$ factors
through $K(S^1)$.
\end{theorem}
\noindent{\em Proof:}
By inspecting the proof of Theorem 6.7, $K(S^1)$ is a minimal simplicial group. The assertion (2) follows
from the following statement.\\
\par
\noindent{\em Statement:} $K(S^1)$ is the quotient simplicial group of $F(S^1)$ modulo the normal subsimplicial group
generated by the boundaries.\\
\par
Let $H$ denote the kernel of the quotient map $p:F(S^1)\rightarrow K(S^1)$ and let $\overline{B}$ denote
the normal subsimplicial group of $F(S^1)$ generated by the boundaries ${\cal B}F(S^1)$.  Notice that
$K(S^1)$ is a minimal simplicial group. Thus $\overline{B}$ is contained in $H$.  Let $Q$ denote
the quotient simplicial group $F(S^1)/\overline{B}$. Then $Q$ is a minimal simplicial group.
By Proposition 2.6, there is a central extension
$$
0\rightarrow K(\pi_{n+1}Q, n+1)\rightarrow P_{n+1}Q\rightarrow P_nQ\rightarrow0,
$$
where $P_nQ$ is the n-th Moore-Postnikov system of $Q$. Notice that $P_1Q=K(\pi_Q,1)=K({\bf Z},1)$.
Thus $\Gamma^{n+2}P_{n+1}Q=1$ by induction on $n$. Notice that $Q_{n+1}\cong (P_{n+1}Q)_{n+1}$.
Thus $\Gamma^{n+2}Q_{n+1}=1$. Now we show that
$H$ is contained in $\overline{B}$ by induction on the dimension starting with $H_1=\overline{B}_1=1$.
Suppose that $H_n\subseteq\overline{B}_n$ with $n>0$. Notice that
$F(S^1)_{n+1}=F(x_0,\cdots,x_n)$ and $K(S^1)_{n+1}=K(x_0,\cdots,x_n)$. Thus $H_{n+1}$ is a normal subgroup of
$F(x_0,\cdots,x_n)$ generated by the commutators
$$
[[x_{i_1},x_{i_2}],\cdots,x_{i_t}],x_{i_t}]
$$
such that $i_p\not=i_q$ for $p<q$. Now consider $W=[[x_{i_1},x_{i_2}],\cdots,x_{i_t}],x_{i_t}]$.
If $t\geq n+1$, then $W\in\Gamma^{n+2}F(x_0,\cdots,x_n)$. Thus $W\in\overline{B}_{n+1}$ since
$\Gamma^{n+2}Q_{n+1}=1$.
\par
If $t<n+1$, then there exists an index $j\in\{0,1,\cdots,n\}-\{i_1,\cdots,i_t\}$. Recall that
$$
s_ix_k=\cases{x_k&$k<i$,\cr x_{k+1}&$k\geq i$,\cr}
$$
for $x_k=s_{n-1}\cdots{\hat s_k}\cdots s_0\sigma\in S^1_n$. Thus
$$
s_j[[x_{i'_1}, x_{i'_2}],\cdots,x_{i'_t}],x_{i'_t}]=[[x_{i_1},x_{i_2}],\cdots,x_{i_t}],x_{i_t}],
$$
where $i'_k=i_k$ if $i_k<j$ and $i'_k=i_k-1$ if $i_k>j$. By induction, 
$[[x_{i'_1}, x_{i'_2}],\cdots,x_{i'_t}],x_{i'_t}]\in\overline{B}$. Thus $W=[[x_{i_1},x_{i_2}],\cdots,x_{i_t}],x_{i_t}]\in\overline{B}$.
The induction is finished and the assertion follows.\\
\par
The simplicial group $K(S^1)$ is homotopy eqivalent to a product of the Eilenberg-Maclane spaces with a different product structure.
\begin{proposition}
$\Omega K(S^1)$ is an abelian simplicial group. Therefore $K(S^1)$ is homotopy equivalent to a product of the Eilenberg-MacLane spaces as a simplicial set.
\end{proposition}

\noindent{\em Proof:} Consider 
$$
d_0:K(x_0,x_1,\cdots,x_n)\rightarrow K(x_0,x_1,\cdots,x_{n-1})
$$
$
d_0(x_0)=1
$
and
$
d_0(x_j)=x_{j-1}.
$ 
Thus  $Kerd_0\cap K_{n+1}(S^1)\cong <x_0>$ is the normal subgroup generated by $x_0$ which is abelian
by Lemma 6.4. The assertion follows.\\
\par
In the end of this section, we give some applications of $K(S^1)$ to minimal simplicial groups.
\begin{proposition}
Let $G$ be a minimal simplicial group such that the abelianlizer $G^{ab}$ is a minimal simplicial group
$K(\pi,1)$ for a cyclic group $\pi$. Then $G$ is homotopy equivalent to a product of Eilenberg-Maclane spaces.
\end{proposition}
\noindent{\em Proof:}
Notice that $G_1=\pi$. Let $x$ be a generator for the cyclic group $\pi$ and let $f_x:S^1\rightarrow G$ 
be a representive map of $x$, i.e, $f_x(\sigma)=x$. Let $g:F(S^1)\rightarrow G$ be the simplicial homomorphism
induced by $f_x$. We need a lemma.
\begin{lemma}
The simplicial homomorphism $g:F(S^1)\rightarrow G$ is simplicial surjection.
\end{lemma}
\noindent{\em Proof:} It suffices to show that the subsimplicial group, denote by $H$, of $G$  generated by $G_1$ is $G$ itself. This is given by induction on
the dimensions starting with $H_1=G_1$. Suppose that $H_{n-1}=G_{n-1}$ with $n>1$. By a result of Condule [see, e.g, Po, Proposition 1, pp.6],
$G_n$ is generated by the degenerate images of lower order Moore chain complex terms and $NG_n$. Thus
$G_n$ is genereted by $NG_n$ and $H_n$ by induction. Notice that $G$ is a minimal simplicial group.
Thus $NG_n={\cal Z}G_n$, the cycles. By Proposition 2.1, ${\cal Z}G_n$ is contained in the center of $G_n$.
Thus $H_n$ is a normal subgroup of $G_n$ and the composite $\phi:\pi_nG\cong{\cal Z}G_n\rightarrow G_n\rightarrow G_n/H_n$ is
an epimorphism. Thus $G_n/H_n$ is an abelian group and so the quotient homomorphism $G_n\rightarrow G_n/H_n$ factors
through $G_n^{ab}$. Notice that $G^{ab}=K(\pi,1)$. Thus $G^{ab}_n\cong G^{ab}_n/{\cal Z}G_n$ and so
$\phi:\pi_nG\rightarrow G_n/H_n$ is trivial. Thus $G_n/H_n$ is trivial and the assertion follows.\\
\par
\noindent{\em Continuation of Proof of Proposition 6.14:}
Notice that $G$ is minimal. The simplicial epimorphism $g:F(S^1)\rightarrow G$ factors through $K(S^1)$ by
Proposition 6.12. By Proposition 6.13, $\Omega K(S^1)$ is an abelian simplicial group. Thus $\Omega G$ is also
an abelian simplicial group. Thus $G$ is homotopy equivalent to a product of Eilenberg-Maclane spaces, which is the assertion.\\
\par
The following counter-example for minimal simplicial groups is due to J. W. Milnor (unpublished).
\begin{proposition}
$\Omega (S^{n+1}[n+1,n+2,n+3])$ does not have a homotopy type of a minimal simplicial group for $n>0$, 
where $S^{n+1}[n+1,n+2,n+3]$ is the 3-stage Postnikov system by taking the first three nontrivial homotopy
groups of $S^{n+1}$.
\end{proposition}
\noindent{\em Proof:} Suppose that $G$ is a minimal simplicial group such that $G\simeq\Omega(S^{n+1}[n+1,n+2,n+3])$.
 Let $f:F(S^1)\rightarrow \Omega^{n-1}G$ be a simplicial
homomorphism such that $f(\sigma)$ is a generator of $(\Omega^{n-1}G)_1\cong G_n\cong {\bf Z}$. Then
$f_*:\pi_j(F(S^1))\rightarrow\pi_j(\Omega^{n-1}G)$ is an isomorphism for $j\leq3$. Notice that
 $\Omega^{n-1}G$ is also a minimal simplicial group. The simplicial homomorphism $f:F(S^1)\rightarrow\Omega^{n-1}G$
factors through $K(S^1)$. Notice that $\pi_3(F(S^1))\cong\pi_3\Omega^{n-1}G)\cong{\bf Z}/2$ and
$\pi_3(K(S^1))\cong Lie(3)\cong{\bf Z}\oplus{\bf Z}$. There is a contradiction and the assertion follows.\\
\par
More examples and counter-examples for minimal simplicial groups will be given in [Wu2]. It was know 
that there are still many counter-examples of two-stage Postnikov systems for minimal simplicial groups [Wu2].

\end{document}